\documentclass[12pt]{amsart}
\usepackage{graphicx}
\usepackage{verbatim}

\newtheorem{theorem}{Theorem}[section]
\newtheorem{lemma}[theorem]{Lemma}
\newtheorem{proposition}[theorem]{Proposition}
\newtheorem{corollary}[theorem]{Corollary}

\theoremstyle{definition}
\newtheorem{example}[theorem]{Example}

\theoremstyle{remark}
\newtheorem{remark}[theorem]{Remark}
\newtheorem*{ack}{Acknowledgments}
\theoremstyle{definition}
 \newtheorem{definition}[theorem]{Definition}

\newcommand{\Tor}{\operatorname{Tor}}
\newcommand{\Int}{\operatorname{Int}}

\newcommand{\lk}{\operatorname{lk}}

\begin{document}

\title{Crosscap Numbers of Two-component Links}
\author{Gengyu Zhang}
\address{Department of Mathematics, Tokyo Institute of Technology, Tokyo 152-8551, Japan}
\email{zhang@math.titech.ac.jp} 
\date{\today}
\thanks{The author is supported by a Research Fellowship of the Japan Society for the Promotion of Science for Young Scientists.
This work is partly supported by the Grant-in-Aid for JSPS Fellows, the Ministry of Education, Culture, Sports, Science and Technology, Japan.
}

\begin{abstract}
We define the crosscap number of a 2-component link as the minimum of the first Betti
 numbers of connected, non-orientable surfaces
 bounding the link. We discuss some properties of the crosscap numbers of 2-component links.
\end{abstract}


\maketitle

\section{Introduction}
Throughout this paper we work in the piecewise linear category, and knots and links we work with are  
embedded in the 3-sphere $S^3$. The crosscap number of a knot $K$ was introduced by Clark \cite{Clark:crosscap} in 1978.
It is defined to be the minimum of the first Betti numbers of non-orientable surfaces bounding $K$. Various notations for 
the crosscap number of a knot have been used in the past research on it, see for example \cite{Clark:crosscap,
Murakami/Yasuhara:crosscap, Teragaito:crosscap}, and in this paper we denote it by $\gamma(K)$.

Clark proved in \cite{Clark:crosscap} the inequality $\gamma(K)\leq 2g(K)+1$ and raised the question whether some knots exist for 
which the equality holds. 
 Murakami and Yasuhara \cite{Murakami/Yasuhara:crosscap} brought a concrete calculation for the knot $7_4$ which 
is the first example known to satisfy the equality above. It has been shown \cite {Hirasawa/Teragaito:crosscap} that
 there exist numerous knots for which the equality holds. 

Given a knot, generally it is hard to determine the crosscap number for it. Clark gave a necessary and sufficient condition 
for the crosscap number 1 knots, which says that a knot has crosscap number 1 if and only if it is a $(2, n)$-cable knot.
Recently the crosscap numbers for several families of knots, such as the torus knots in 
\cite{Teragaito:crosscap}, the 2-bridge knots in \cite{Hirasawa/Teragaito:crosscap}, and the pretzel knots in 
\cite{Ichihara/Mizushima:crosscap}, have been 
determined.

In this paper we define the crosscap number for two-component links and discuss some properties of it. By following the 
technique used in \cite{Murakami/Yasuhara:crosscap}, we calculate the crosscap number of the two-component link $6_3^2$ as
an example. Here we use the notation of Rolfsen \cite{Rolfsen:1990} to denote a link in his link table.

\section{definitions}

The crosscap number of a knot was first introduced by Clark \cite{Clark:crosscap} in 1978. 
\begin{definition}[\cite{Clark:crosscap}]
The \textit{crosscap number} $\gamma(K)$ of a knot $K$ is the minimal number of the first Betti numbers of all the connected, non-orientable
surfaces bounding $K$. The crosscap number of the unknot $U$ is defined to be 1.
\end{definition}

Note that in this paper we define the crosscap number of the unknot to be 1, 
instead of 0 defined by Clark \cite{Clark:crosscap}.

Clark also gave an upper bound for the crosscap number of a knot in terms of its genus.
\begin{proposition}[\cite{Clark:crosscap}]
Let $K$ be a knot, and $g(K)$ denote the genus of $K$. Then
\begin{equation}
\gamma(K) \leq 2g(K)+1.
\end{equation}
\end{proposition}

Beginning with the knot $7_4$ proved by Murakami and Yasuhara in \cite{Murakami/Yasuhara:crosscap}, it has been 
shown in \cite{Hirasawa/Teragaito:crosscap} that numerous knots are suited for the equality in (1).

The crosscap number of a two-component link is defined similarly to that of a knot. 
\begin{definition}
The \textit{crosscap number} $\gamma(L)$ of a two-component link $L$ is the minimum of the first Betti numbers of connected, 
non-orientable surfaces bounding the link, i.e. we have
 \begin{multline*}
  \gamma(L):=\min \{ \beta_1(F)
 \\
  \mid \text{$F$ is a connected non-orientable surface bounding $L$}\}.
  \end{multline*}
\end{definition}

\begin{figure}
   \includegraphics [scale=.4] {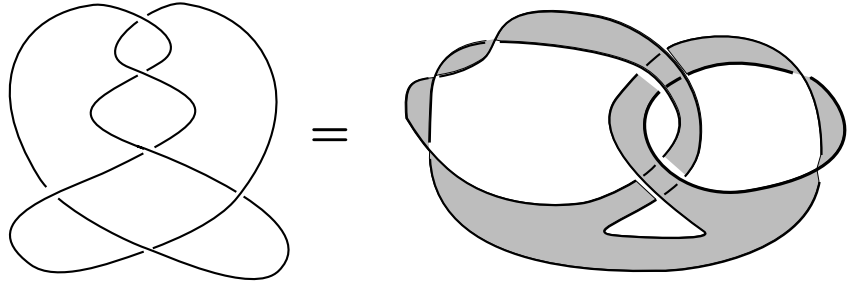}
   \vspace {-7.3cm}
   \caption{Non-orientable surface bounding $6_2^2$ with $\beta_1=2$}
   \label{fig:betti_number_2_622}
\end{figure}
It is not hard to see that for a 2-component link $L$, its crosscap number is at least 2, i.e. $\gamma(L)\geq 2$. 
This is because the projective plane $\mathbf{R}P^2$ is a closed surface with minimum first Betti number, and 
the first Betti number of the surface obtained by cutting two disks off $\mathbf{R}P^2$ is 2.   

Let $6_2^2$ be the two-component link illustrated in Figure~\ref{fig:betti_number_2_622}.
\begin{example} 
We have $\gamma(6_2^2)=2$ as shown in the picture.
\end{example}

\section{Behavior of crosscap numbers under split union}

In this section we will see how the crosscap number of a split union of two knots can 
be evaluated by the crosscap numbers of the knots. 

\begin{definition}[See \cite{Lickorish:knot}]
A link $L=L_1\cup L_2$ is called \textit{splittable} if there exists some 
2-sphere $S^2$ embedded in $S^3$ bounding two 3-balls $B_1$ and $B_2$ with $B_1\cup B_2=S^3$ such that $L_i$
lies in the interior of $B_i$ $(i=1,2)$. A \textit {split union} of two knots $K_1$ and $K_2$, denoted by $K_1\circ K_2$,
 is the splittable link with the 2-sphere bounding two 3-balls such that $K_i$ lies in the interior of $B_i$ ($i=1,2$). 
\end{definition}

Firstly, we consider the crosscap number of the split union of any two knots.
In actuality, for the split union of any two knots, its crosscap number can be known in terms of genera and crosscap numbers 
of the knots. Namely, we have the following equality.

\begin{theorem} \label{thm:min}
Let $K_1$ and $K_2$ be any two knots. We have
      \begin{multline*}
      \gamma(K_1\circ K_2)
        =\min\{\gamma(K_1)+2g(K_2)+1, \\\gamma(K_2)+2g(K_1)+1, \gamma(K_1)+\gamma(K_2)+1\}. 
      \end{multline*}
\end{theorem}
\begin{proof}
Let $G'$ be a connected, non-orientable surface which gives the minimum first Betti number
 for the split union $K_1\circ K_2$. 
 There exists a 2-sphere $S$ separating $K_1$ and $K_2$ with $S\cap G'\neq \emptyset$. 
 Take an innermost cirle $\alpha$ of  $S\cap G'$ in $S$ without bounding any disk in $G'$.  
(We can always make it 
because we can remove the circle $\alpha$ without changing the first Betti number by doing surgery on the surface $G'$ 
if it bounds a disk in $G'$.) 

 If $\alpha$ is
 a non-separating curve in $G'$, then we do surgery along $\alpha$ on the surface $G'$, from which a new surface with smaller
 first Betti number, still bounding $K_1\circ K_2$, connected, will appear. This contradicts the assumption that 
 the surface $G'$ realizes the crosscap number for $K_1\circ K_2$. 
 
 Hence we may assume that $\alpha$ is a separating curve which
 separates the surface $G'$ into two surfaces, say $G_1'$ and $G_2'$. 
 Each of them bounds $K_1$ or $K_2$ separately; otherwise we can get a connected, 
non-orientable surface bounding $K_1\circ K_2$ with smaller first Betti number by doing surgery along $\alpha$ on $G'$,
which contradicts the minimality of the crosscap number. 
 Then at least one of the two surfaces $G_1'$ and $G_2'$ is non-orientable. 
 In all there are three possibilities:
 $G_1'$ orientable and $G_2'$ non-orientable, $G_1'$ non-orientable and $G_2'$ orientable, 
 or both non-orientable.
 Therefore we get the inequality $\min\{\gamma(K_1)+2g(K_2)+1, \gamma(K_2)+2g(K_1)+1, \gamma(K_1)+\gamma(K_2)+1\}\leq \gamma(K_1\circ K_2)$.

On the other hand, we may assume that the knot $K_i$ bounds an orientable surface $S_i$ with genus $g(K_i)$ and 
a non-orientable surface $G_i$ with first Betti number $\gamma(K_i)$ ($i=1,2$). Then three non-orientable surfaces bounding
$K_1 \circ K_2$ will be produced if we connect the surfaces $S_1$ and $G_2$, $G_1$ and $S_2$, and $G_1$ and $G_2$ by tubes.
The first Betti numbers of these three surfaces are $\gamma(K_1)+2g(K_2)+1$, $\gamma(K_2)+2g(K_1)+1$ and $\gamma(K_1)+\gamma(K_2)+1$
respectively. This gives us an upper bound of the crosscap number of the split union,
i.e. $\gamma(K_1\circ K_2)\leq \min\{\gamma(K_1)+2g(K_2)+1, \gamma(K_2)+2g(K_1)+1, \gamma(K_1)+\gamma(K_2)+1\} $.

Then we have the proof of the equality.
\end{proof}

\begin{corollary}\label{thm:inequality}
 Let $K_1$ and $K_2$ be any two knots. 
 Then the following inequalities hold:
   \begin{align}
    \gamma(K_1\circ K_2)\leq \gamma(K_1)+\gamma(K_2)+1,\\ \notag
    \gamma(K_1\circ K_2)\leq \gamma(K_1)+2g(K_2)+1,\\ \notag
    \gamma(K_1\circ K_2)\leq \gamma(K_2)+2g(K_1)+1.
    \end{align}
\end{corollary}

Using Clark's inequality (1), we have
\begin{corollary}\label{thm:sum}
Let $K_1$ and $K_2$ be any two knots. 
 Then $\gamma(K_1 \circ K_2)=\gamma(K_1)+\gamma(K_2)+1$ if and only if $\gamma(K_i)< 2g(K_i)+1$, $i=1$ and $2$.
\end{corollary}

This corollary is equivalent to the following.

\begin{corollary}
$\gamma(K_1 \circ K_2)=\gamma(K_1)+\gamma(K_2)$ if and only if $\gamma(K_i)= 2g(K_i)+1$, $i=1$ or $2$.
\end{corollary}

Note that when $\gamma(K_1)<2g(K_1)+1$ and $\gamma(K_2)=2g(K_2)+1$, then $\gamma(K_1\circ K_2)= \gamma(K_1)+2g(K_2)+1$;
          when $\gamma(K_1)=2g(K_1)+1$ and $\gamma(K_2)<2g(K_2)+1$, then $\gamma(K_1\circ K_2)= \gamma(K_2)+2g(K_1)+1$. 
In fact, the crosscap number of $K_1\circ K_2$ is exactly equal to $\gamma(K_1)+\gamma(K_2)$ in both cases.

If we apply the argument above to the case when $K_2=U$, we have the following corollary.
\begin{corollary}
 Let $U$ denote the unknot and $K$ be any knot. Then we have
   \begin{equation*}
     \gamma(K\circ U)=\gamma(K)+1.
   \end{equation*}
\end{corollary}

Next, by applying the homology theory, we discuss a little more the examples for which the equality in (2) holds. 
\begin{proposition} \label{prop:order}
Let $D(L)$ denote the double branched cover of $S^3$ branched along the two-component link $L$.  
Then the minimum number of generators for $H_1(D(L);\mathbf{Z})$ has the crosscap number $\gamma(L)$ as an upper bound.
\end{proposition}

\begin{proof}
Let $F$ denote a non-orientable surface which has the minimum first Betti number, $\gamma(L)$, bounding $L$. Then corresponding to this surface, 
there is a $\gamma(L)\times \gamma(L)$ Goeritz matrix built in the way of Gordon and Litherland \cite[\S 2]{Gordon/Litherland:signature}.
Then this Goeritz
matrix becomes a relation matrix for $H_1(D(L);\mathbf{Z})$, see Appendix~\ref{sec:a}, from which the result follows. 
\end{proof}

\begin{figure}
\includegraphics [scale=.4] {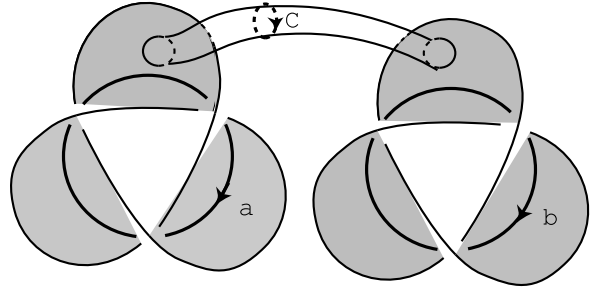}
\vspace {-5.5cm}
 \caption{Tubular connection}
\label{fig:tubular_connection}
\end{figure}

\begin{remark} \label{thm:trefoils}
It is known that $g(3_1)=1$ and $\gamma(3_1)=1$. Thus $\gamma(3_1)< 2g(3_1)+1$, and by Corollary~\ref{thm:sum} we have
$\gamma(3_1\circ 3_1)=3=\gamma(3_1)+\gamma(3_1)+1$. Now we prove this in terms of homology theory.

By connecting the two non-orientable surfaces bounding trefoils with a tube, a non-orientable surface $F$
  bounding $3_1\circ 3_1$ can be built. see Figure~\ref{fig:tubular_connection}.
The Georitz matrix corresponding to this surface with indicated generators \{$a,b,c$\} is as follows:
\[ \left(
\begin{array}{lcr}
3&0&0\\0&3&0\\0&0&0 
\end{array} 
\right). \]

Then the first homology group of the double cover of $S^3$ branched over $3_1\circ 3_1$,  
can be known as $H_1(D(3_1\circ 3_1);\mathbf {Z})=\mathbf{Z}/3\mathbf{Z}\oplus \mathbf{Z}/3\mathbf{Z}\oplus \mathbf{Z}$. By the fundamental theorem of abelian groups, $H_1(D(3_1\circ 3_1);\mathbf{Z})$ cannot be presented by a $2\times 2$ matrix.

So the crosscap number of the split union of two trefoil knots cannot be 2 by Proposition~\ref{prop:order}. 
Then we have $\gamma(3_1\circ 3_1)=3=\gamma(3_1)+\gamma(3_1)+1$.
\end{remark}

\section{upper bounds of crosscap numbers of two-component links}

Denote by $\overrightarrow{L}$ and $\overleftarrow{L}$ the two different relative orientations 
for a 2-component link $L$. Let $g(\overrightarrow{L})$ and $g(\overleftarrow{L})$ denote the genera of $L$
under these two different orientations.

\begin{theorem}
With the notations above, we have
\begin{equation}
\gamma(L)\leq 2\min(g(\overrightarrow{L}),g(\overleftarrow{L}))+2.
\end{equation}
\end{theorem}

\begin{figure}
\includegraphics [scale=.33] {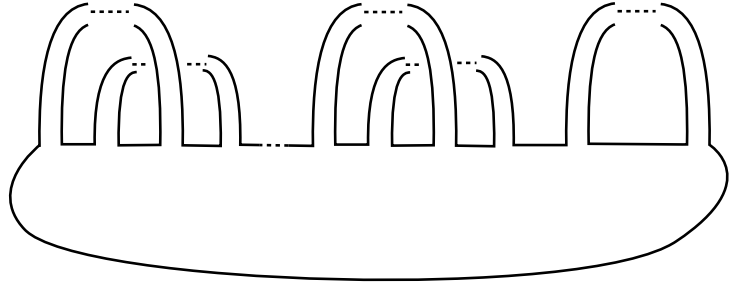}
\vspace {-1.8cm}
 \caption{Standard form of a Seifert surface bounding a 2-component link}
\label{fig:standard_form}
\end{figure}
\begin{figure}
\includegraphics [scale=.33] {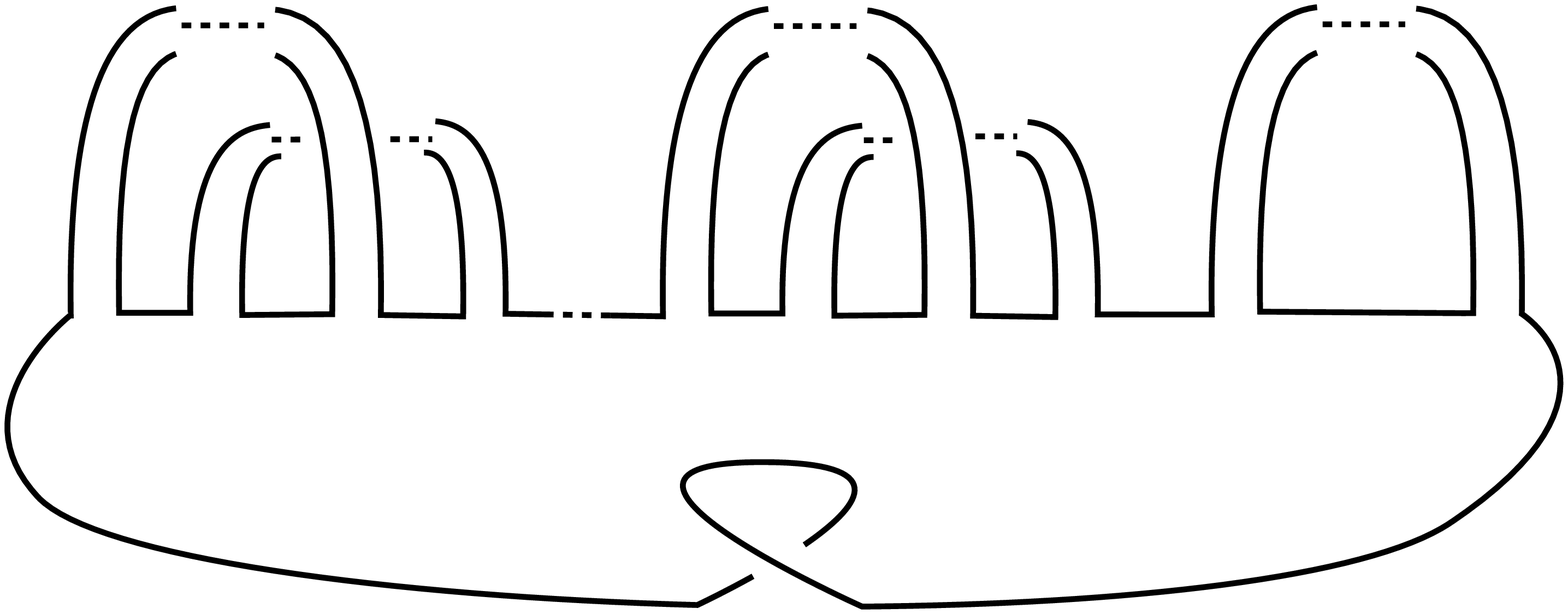}
\vspace {-5.4cm}
 \caption{Surface after adding a half-twisted band}
\label{fig:nonorientable}
\end{figure}

\begin{proof}
Denote $\min(g(\overrightarrow{L}),g(\overleftarrow{L}))$ by $g$. Then there exists an orientable Seifert surface $F$, whose standard form is illustrated in 
Figure~\ref{fig:standard_form}, with genus $g$ bounding the 2-component link $L$, so the first Betti number of this surface becomes $2g+1$.
By adding a half twist to the surface, 
we obtain a non-orientable surface from $F$ as illustrated in Figure~\ref{fig:nonorientable}, 
whose first Betti number is $2g+2$. By the definition of the crosscap number, we know that the inequality holds.
\end{proof}

Does there exist any $2$-component link for which the equality in {\rm (3)} holds? 

\begin{figure}
\includegraphics [scale=.25] {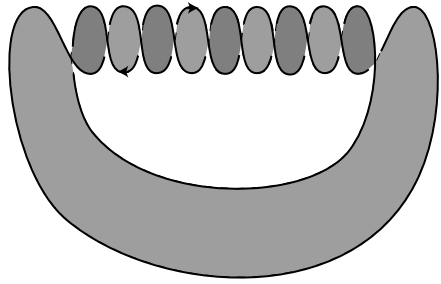}
\vspace{-4.5cm}
\caption{Torus link $T(2,10)$}
\label{fig:torus2_10link}
\end{figure}

There exists an infinite sequence 
of 2-component links $T(2,2n)$ ($n\in \mathbf{Z}$) for which the equality in {\rm(3)} holds, where $T(p,q)$ denotes a torus 
knot or link. Take an example of torus link $T(2, 10)$ 
as illustrated in Figure~\ref{fig:torus2_10link}. 
It bounds a genus $0$ orientable Seifert surface, which gives us an upper bound $2$ for the crosscap number of 
the link by using the inequality {\rm(3)}. 
Therefore we have the fact that the crosscap number of the torus link $T(2,2n)$ is 2, for which the equality in {\rm(3)} holds. 

Let $n(L)$ denote the minimum crossing number of a link $L$. We can give another upper bound
for the crosscap number of a two-component link in terms of $n(L)$. For a knot,
 Murakami and Yasuhara in \cite{Murakami/Yasuhara:crosscap} proved the following proposition.
\begin{proposition}[\cite{Murakami/Yasuhara:crosscap}] \label{thm:upper}
Let $n(K)$ denote the crossing number of a knot $K$. Then
\begin{equation}
\gamma(K) \leq \dfrac{\lfloor n(K) \rfloor}{2}, 
\end{equation}
where $\lfloor x \rfloor$ denotes the greatest integer that does not exceed $x$.
\end{proposition}

 For a two-component link, we have the following result.
\begin{proposition}\label{thm:crossing}
Let $L$ be a two-component link, excluding the unlink. Then we have 
\begin{equation*}
\gamma(L) \leq \dfrac{\lfloor n(L) \rfloor}{2}+1,
\end{equation*}
where $\lfloor x \rfloor$ denotes the greatest integer that does not exceed $x$.
\end{proposition}

Note that for the Hopf link, the equality in Proposition~\ref{thm:crossing} holds.

\begin{proof}
If $L$ is a splittable link $K_1 \circ K_2$, then from Proposition~\ref{thm:upper} we have
$\gamma(K_i) \leq \dfrac{\lfloor n(K_i) \rfloor}{2}$ for $i=1,2$, and from Corollary~\ref{thm:inequality}
we have $\gamma(L)\leq \gamma(K_1)+\gamma(K_2)+1$. Then the inequality
$\gamma(L)\leq \dfrac{\lfloor n(K_1) \rfloor}{2}+\dfrac{\lfloor n(K_2) \rfloor}{2}+1$ follows. Hence 
the theorem holds for splittable links due to the fact that $\dfrac{\lfloor n(K_1) \rfloor}{2}+\dfrac{\lfloor n(K_2) \rfloor}{2}\leq \dfrac{\lfloor n(L) \rfloor}{2}$.

Now let $D$ be a link diagram of non-splittable link $L$ with the minimum crossing number $n(L)$. Then we have $n(L)+2$ regions of $S^2$ divided by 
the link diagram. Color these regions black and white in a checkerboard way. Since $L$ is a non-splittable link, 
all the regions of the same color can be connected to each other by half-twisted bands at the crossings. 
Both the white and black surfaces can be orientable,
then $\gamma(L)$ should be less than or equal to the first Betti number of these surfaces plus one, where the ``1''
indicates the added first Betti number by adding a half-twisted band in that case.
Denote the numbers of 
black and white regions by $n(b)$ and $n(w)$ respectively. Note that the number of the edges is twice that of vertices.
Then by using a relation between Euler characteristic and the first 
Betti number together with Euler's formula, we have 
\begin{equation*}
\gamma(L) \leq 2+n(L)-\max\{n(b), n(w)\}.
\end{equation*}
It is not hard to know that $\max\{n(b), n(w)\} \geq \frac{1}{2}(n(L)+2)$ in the case that $n(L)$ is even, and
$\max\{n(b), n(w)\} \geq \frac{1}{2}(n(L)+3)$ in the case that $n(L)$ is odd. Therefore the result follows. 
\end{proof}

\section{An example of calculation}

In this section, we will calculate the crosscap number of the two-component link $6_3^2$ as an example.

\begin{lemma}\label{thm:2x2}
Let $L$ be a two-component link $K_1 \cup K_2$. Assume that it bounds a connected, non-orientable surface $F$ with the first 
Betti number $2$. Then we can choose a generator system for $H_1(F;\mathbf{Z})$ such that the Goeritz matrix $G_F(L)$ corresponding to this system is of the following form:

\begin{equation}    
G_F(L)=\begin{pmatrix}
2n+1&2k\\2k&2m \end{pmatrix}
\end{equation},
where $k, m, n$ are integers.
\end{lemma}

\begin{figure}
 \includegraphics [scale=.4] {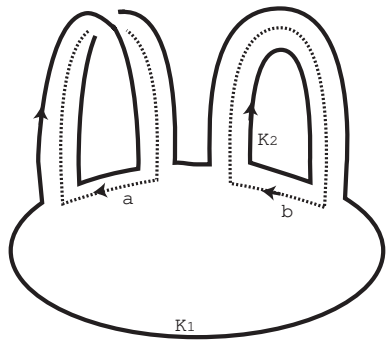}
 \caption{A crosscap number two surface}
\label{fig:generator_system}
\end{figure}

\begin{figure}
\includegraphics [scale=.3] {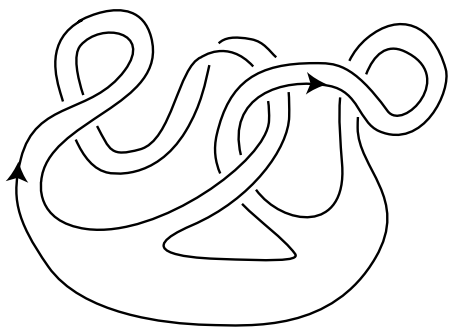}
\vspace {-2.8cm}
\caption{An example of a non-orientable surface with $\beta_1=2$}
\label{fig:example_2_surface}
\end{figure}

\begin{proof}
We may assume that $F$ is a disk with a non-orientable band and an orientable one as indicated in Figure~\ref{fig:generator_system},
where each band may be knotted and linked with each other; see for example Figure~\ref{fig:example_2_surface}. Choose a generator system \{$a, b$\} as in Figure~\ref{fig:generator_system} and orient the 
two components $K_1$ and $K_2$ so that the two boundaries of each band have the same orientations as that of the 1-cycle passing through it.

Then the Goeritz matrix of the surface $F$, refer to \cite[\S 2]{Gordon/Litherland:signature}, corresponding to this generator system is
\begin{equation*}
\begin{pmatrix} 
\lk(a, \tau a)&\lk(a,\tau b) \\
\lk(b,\tau a)&\lk(b,\tau b) 
        \end{pmatrix}
=\begin{pmatrix}
2\delta  (A, A)+1&2\delta (A, B)\\
2\delta (B, A)&2\delta (B, B)
\end{pmatrix},
\end{equation*}
where $\tau (x)$ denotes the orientation double cover of a cycle $x$. 

\begin{figure}
\includegraphics [scale=.35] {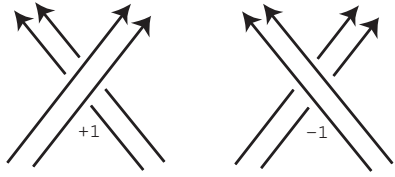}
\caption{Crossing types of bands}
\label{fig:band_crossing_type}
\end{figure}

Here $A$ and $B$ are the bands which $a$ and $b$ pass through respectively, and $\delta (X, Y)$ is the sum of signs of crossings 
of bands $X$ and $Y$ with signs determined as in Figure~\ref{fig:band_crossing_type}. Note that $\lk(a, \tau b)=\lk(b, \tau a)=2\lk(a, b)$, 
and that both of them are even. This finishes the proof.
\end{proof}

\begin{lemma} \label{thm:linking_number_Euler_number}
With the orientations as above in the proof of Lemma~\ref{thm:2x2}, the linking number between $K_1$ and $K_2$ is $m+2k$ 
and the modified Euler number, see Appendix~\ref{sec:a}, of the surface $F$ is $(-2)$ times the sum of all 
the elements in the Goeritz matrix $G_F(L)$. Namely,
\begin{equation*}
\lk(K_1, K_2)=m+2k, \quad \overline{e}(F)=-2(2n+1+2k+2k+2m).
\end{equation*}
\end{lemma}

\begin{figure}
  \includegraphics [scale=.49] {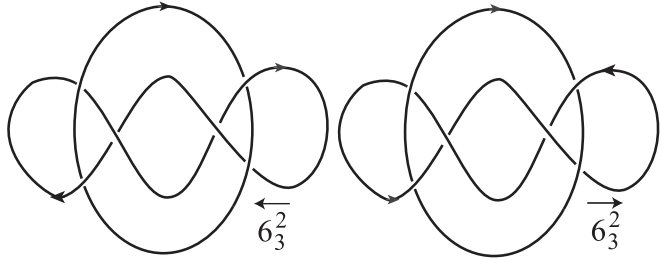}
 \caption{2-component link $6_3^2$ with two orientations}
\label{fig:6322componentlink}
\end{figure}

\begin{proof}
According to the proof of the above theorem, we have $m=\delta(B, B)$ and $k=\lk(a, b)$. It is obvious that the calculation of
the linking number between $K_1$ and $K_2$ includes these two parts, which gives us $\lk(K_1, K_2)=\delta(B, B)+2\lk(a,b)=m+2k$.

We also have the fact that the modified normal Euler number of $F$ with the orientations above is 
$-\{4\delta (A, A)+2+4\delta (A, B)+4\delta (B, A)+4\delta (B, B)\}$. 
Then since $m=\delta(B, B)$, $k=\delta(A, B)$ and $n=\delta(A,A)$, we obtain the second equality.
\end{proof}

\begin{figure}
  \includegraphics [scale=.35] {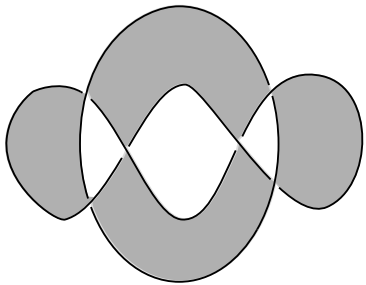}
  \vspace{-6.3cm}
 \caption{Non-orientable surface with $\beta_1=3$}
\label{fig:bettinumber3632}
\end{figure}
We will prove the following conclusion.

\begin{theorem}
Let $6_3^2$ denote the two-component link as illustrated in Figure~\ref{fig:6322componentlink}, forgetting the orientations.
 We have 
$\gamma(6_3^2)=3$.
\end{theorem}
\begin{proof}
It is clear that the surface bounding the two-compnent $6_3^2$ 
illustrated in Figure~\ref{fig:bettinumber3632} is a non-orientable one with first Betti number 3. 

Suppose that $\gamma(6_3^2) \leq 2$. Then there exists a non-orientable surface $F$ bounding the link
 with $H_1(F; \mathbf{Z})=\mathbf {Z}\oplus \mathbf{Z}$. The Goeritz matrix $G_F(6_3^2)$ associated with $F$ should be
a 2 by 2 matrix, determined by a generator system of $H_1(F; \mathbf{Z})$. A different choice of basis of 
 $H_1(F; \mathbf{Z})$ gives another matrix $J$ such that $J=P^TGP$ where $P$ is an integral unimodular matrix.
 The integral congruent class which $G_F(6_3^2)$ belongs to does not change. 

Since the link $6_3^2$ is the two-bridge link 
$S(12,5)$, then the double branched cover $D(6_3^2)$ of $S^3$ branched over $6_3^2$ is the lens space $L(12, 5)$ with
$H_1(D(6_3^2); \mathbf{Z})=\mathbf{Z}/12\mathbf{Z}$ and the linking form $\lambda (g,g)=\pm \frac{5}{12}$ for some properly 
chosen generator of $H_1(D(6_3^2); \mathbf{Z})$.
The determinant of the Goeritz matrix $G_F(6_3^2)$ is known to be $\pm 12$ since the absolute value of the determinant of 
a Goeritz matrix of a link is equal to the order of $H_1(D(L); \mathbf{Z})$. 

By applying an elementary theorem of integral binary quadratic forms (see, for example \cite{Niven:numbers}), we enumerate all the congruent classes of 2 by 2
integral matrices with discriminant $\pm 48$. Then the result is as follows:
\begin{gather*}
X_1=\begin{pmatrix}1&0\\0&12\end{pmatrix}, X_2=\begin{pmatrix} 2&0\\0&6 \end{pmatrix}, X_3=\begin{pmatrix} 3&0\\0&4 \end{pmatrix},
X_4=\begin{pmatrix}4&2\\2&4 \end{pmatrix},\\ X_5=\begin{pmatrix} 1&0\\0&-12 \end{pmatrix}, X_6=\begin{pmatrix} -1&0\\0&12 \end{pmatrix},
X_7=\begin{pmatrix} 2&0\\0&-6 \end{pmatrix}, \\X_8=\begin{pmatrix} -2&0\\0&6 \end{pmatrix}, X_9=\begin{pmatrix} 3&0\\0&-4 \end{pmatrix},
X_{10}=\begin{pmatrix} -3&0\\0&4 \end{pmatrix}.
\end{gather*} 

It is obvious that the matrices $X_2, X_4, X_7$ and $X_8$ cannot present the cyclic group $\mathbf{Z}/12\mathbf{Z}$,
 and therefore cannot be the relation matrices of $H_1(D(6_3^2); \mathbf{Z})$. On the other hand, only the matrix $X_3$ 
presents the linking form $\pm \frac{5}{12}$ for the link $6_3^2$, see Appendix~\ref{sec:b}.

If $X_3$ were a Georitz matrix for 
the link $6_3^2$, then there should exist some integral unimodular matrix $Q$ such that $Q^TG_F(L)Q=X_3$.
The inverse matrix of $Q$ is also integral unimodular and we denote it by $P$. 
Changing the basis of $H_1(D(6_3^2); \mathbf{Z})$ by using $P:=\begin{pmatrix} r&u\\s&v \end{pmatrix}$, 
the Goeritz matrix $G_F(L)$ is of the form:
\begin{equation*}
P^T\begin{pmatrix}3&0\\0&4\end{pmatrix}P=\begin{pmatrix}3r^2+4s^2&3ru+4sv\\3ru+4sv&3u^2+4v^2\end{pmatrix}.
\end{equation*}

Recall that the formula by Gordon and Litherland \cite{Gordon/Litherland:signature}, also see Appendix~\ref{sec:a}, relates the 
signature of a link with the signature of the Goeritz matrix, 
$\sigma (\overline L)=\sigma(G_F) +\dfrac{1}{2}\overline e(F)$. The signature of a link is defined as the signature of
the symmetrized Seifert matrix, the difference between the number of positive eigenvalues and negative ones of the matrix.

For some appropriate Seifert surfaces, the corresponding Seifert matrices of the two-component link 
under two different orientations, $\overleftarrow {6_3^2}$ and $\overrightarrow {6_3^2}$, 
see Figure~\ref{fig:6322componentlink},
are as follows:
\begin{equation*}
V(\overleftarrow {6_3^2})=
\begin{pmatrix}
2&1&0\\0&1&0\\-1&0&1
\end{pmatrix}, 
V(\overrightarrow {6_3^2})=
\begin{pmatrix}
-1&-1&0\\0&1&0\\0&1&-1
\end{pmatrix}. 
\end{equation*}
Then we obtain the signatures of $\overleftarrow {6_3^2}$ and $\overrightarrow {6_3^2}$ as 3 and $-1$ respectively.
We will see that for either orientation, there exists no solution for which the formula by Gordon and Litherland holds. 

According to Lemma~\ref{thm:linking_number_Euler_number}, we see that the modified normal Euler number $\overline e(F)$ is $-2[3(r+u)^2+4(s+v)^2]$. 
Now since $\sigma (X_3)=2$ we have
\begin{equation}
\sigma (\overleftarrow {6_3^2})=3=2-[3(r+u)^2+4(s+v)^2],
\end{equation}
or 
\begin{equation}
\sigma (\overrightarrow {6_3^2})=-1=2-[3(r+u)^2+4(s+v)^2].
\end{equation}

Since a negative number cannot be equal to the sum of perfect squares, 
it is not hard to see that there exists no integral solution of $r, s, u, v$ to the equality (6). 

For $\overrightarrow {6_3^2}$, we establish another equality in terms of the linking number 
by Lemma~\ref{thm:linking_number_Euler_number}. We have
\begin{equation}
\lk(\overrightarrow {6_3^2})=2=\dfrac{3u^2+4v^2}{2}+3ru+4sv.
\end{equation}
Then we combine the equalities (7) and (8) to get the following system of equations:
\begin{equation*}
\begin{cases}
3(r+u)^2+4(s+v)^2=3,\\
3u^2+4v^2+6ru+8sv-4=0.
\end{cases}
\end{equation*}

Because the values of $r, u, s$ and $v$ are taken in $\mathbf{Z}$, from the first equation we have $r+u=\pm 1$ and $s+v=0$.
Transform the second equation to the following:
\begin{equation*}
3(u+r)^2-3r^2+4(s+v)^2-4s^2-4=0.
\end{equation*}
Putting the values of $u+r$ and $s+v$ into the transformed equation, we have $3-3r^2-4s^2-4=0$, i.e., $3r^2+4s^2=-1$, for which there exists no integral solution of $r$ and $s$.

So it turns out that there exist no integral solutions for this system of equations, which contradicts our assumption.
Namely, there exists no connected, non-orientable surface with first Betti number 2. Therefore, we have $\gamma(6_3^2)=3$.
\end{proof}

\begin{ack}
The author would like to express her deep thanks to Professor Hitoshi Murakami for helpful discussions. She also would like to thank
the Ministry of Education, Culture, Sports, Science and Technology of Japan for offering her MEXT scholarship when the main part 
of the work was carried out. 
\end{ack}

\appendix
\section{A formula of signature}\label{sec:a}

The definition of the signature of knots and links was developed by Trotter \cite{Trotter:homology} and Murasugi \cite{Murasugi:signature}. 
For an oriented link in $S^3$ with Seifert matrix $V$, define the \textit{signature} of $L$ to be $\sigma (L)=\sigma (V+V^T)$ where 
$\sigma (M) $ is the difference between the number of the positive eigenvalues and that of the negative eigenvalues of a symmetric matrix $M$.
 Note that the signature of a link is up to the relative orientations for each component.
See, for example \cite{Rolfsen:1990} or \cite{Lickorish:knot}.

Gordon and Litherland \cite {Gordon/Litherland:signature} has shown how to define a quadratic form related with Goeritz matrix by using any spanning surface, and related the signature
 of this form to the signature of a link. We recall their formula here for readers' information.
\begin{lemma}[\cite {Gordon/Litherland:signature}\label{thm:signature}]
Let $F$ be any surface bounding an unoriented link $L$, and let $\overline L$ denote the link $L$ together with some orientation on each component of it.
 Then the signature $\sigma (\overline L)$ can be calculated out of the Goeritz matrix $G_F$ and the modified normal Euler number $\overline e(F)$, namely we have
\begin{equation}
\sigma (\overline L)=\sigma(G_F) +\dfrac{1}{2}\overline e(F).
\end{equation} 
\end{lemma}

\begin{definition}[\cite {Gordon/Litherland:signature}]
  Let $L=K_1 \cup K_2$ be a $2$-component link with a connected surface $F$ bounding it. 
 We define $\overline e(F)$ to be $ -\sum_{i,j=1}^2 \lk(K_i, K_j^\prime)$  where $K_i^\prime$ is the intersection
 of $F$ and the boundary of the regular neighborhood of $K_i$ in $S^3$, i.e. 
$K_i^\prime =F\cap \partial N(K_i)$, with orientation parallel to that of $K_i$.
We name it the \textit{modified normal Euler number}.
\end{definition}

\section{Linking form of a link} \label{sec:b}
Let $M$ denote a closed, oriented 3-dimensional manifold. Denote the torsion part of $H_1(M;\mathbf{Z})$ by $\Tor H_1(M;\mathbf{Z})$. A chain complex of $M$ is as follows:
\begin{equation*}
\dots\longrightarrow C_2(M) \stackrel{\partial_2} {\longrightarrow} C_1(M) \stackrel{\partial_1}{\longrightarrow} 0,
\end{equation*}
where $C_i(M)$ is the abelian group generated by $i$-simplices of $M$, and $\partial_i$ is the boundary homomorphism ($i= 1, 2, \dots$).
Suppose that $x, y \in \Tor H_1(M;\mathbf{Z})$ are represented by 1-cycles $a$ and $b$ respectively. There exists $n \in \mathbf{Z}$ such that
 $[nb]$ is homologous to zero in $H_1(M;\mathbf{Z})$ and $nb$ forms the boundary of some 2-chain, say $\Delta $ in $C_2(M)$, i.e. $\partial_2(\Delta)=nb$. 
Define a bilinear form $\lambda: \Tor H_1(M;\mathbf{Z})\times \Tor H_1(M;\mathbf{Z})\longrightarrow \mathbf{Q}/\mathbf{Z}$ as follows:
\begin{equation*}
(x, y) \stackrel{\lambda}{\longmapsto} \Int(a, \Delta)/n,
\end{equation*}
where $\Int$ denotes the intersection number between a 1-cycle and a 2-chain. This bilinear form $\lambda$ is called the \textit{linking form} on the 3-manifold $M$.

Then the linking form defined on the double branched cover $D(L)$ of $L$ 
   \begin{equation*}
     \lambda: \Tor H_1(D(L);\mathbf{Z})\times \Tor H_1(D(L);\mathbf{Z})\longrightarrow \mathbf{Q}/\mathbf{Z} 
   \end{equation*}
is called \textit{the linking form of the link} $L$. A Goeritz matrix $G$ of the link $L$ is a relation matrix for $H_1(D(L);\mathbf{Z})$, 
and the first homology group of $D(L)$ and the linking form on $D(L)$, $(H_1(D(L);\mathbf{Z}), \lambda)$, 
can be calculated out of the Goeritz matrix. Precisely speaking, $H_1(D(L);\mathbf{Z})=\mathbf{Z}^n/\operatorname{Im}(G)$ 
where $n$ is the size of the Goeritz matrix $G$, and the linking form is given by $\pm G^{-1}$ if $H_1(D(L);\mathbf{Z})$ is finite
, i.e.  
$\lambda(g_i, g_j)=\pm (G^{-1})_{ij}\pmod 1$ for the generators $g_i$ and $g_j$ of $H_1(D(L);\mathbf{Z})$, 
where \{$g_1, g_2, \cdots, g_n$ \} is a generator system of $H_1(M;\mathbf{Z})$ corresponding to the presentation of 
$\mathbf{Z}^n/\operatorname{Im}(G)$.
Here the sign $\pm$ depends on the orientation.


\bibliography{mrabbrev,gengyu}
\bibliographystyle{amsplain}

\end{document}